\documentclass{article}

\usepackage{proof}
\usepackage{latexsym}

\newtheorem{theorem}{Theorem}

\newtheorem{lemma}[theorem]{Lemma}

\title{Refinement of proof-theoretic analysis of the lpo due to Buchholz
}
\author{Toshiyasu Arai
\\
Graduate School of Mathematical Sciences
\\
University of Tokyo
\\
3-8-1 Komaba, Meguro-ku,
Tokyo 153-8914, JAPAN
\\
tosarai@ms.u-tokyo.ac.jp
}
\date{}
\begin{document}
\maketitle

\begin{abstract}
We give a refinement of proof-theoretic analysis of the lpo (lexicographic path order)
due to W. Buchholz\cite{Buchholz}.
This note was written in Feb.\,5, 2015 when G. Moser visited Japan.
\end{abstract}

For a binary relation $<$ on a set $T$, $W(<)$ denotes the well-founded part of $T$ with respect to $<$.

Let $\mathcal{F}$ be a non-empty and finite set of function symbols.
Suppose that a precedence (an irreflexive and transitive relation) $<$ on $\mathcal{F}$ is given.
Each $f$ has a fixed arity $ar(f)\in\omega$.
Assume that there is a function symbol $f\in\mathcal{F}$ such that $ar(f)>1$.
When writing $f(a)$, we tacitly assume that $a\in T(\mathcal{F})^{ar(f)}$.
If $a=(a_{0},\ldots,a_{n-1})$, then
$t\in a$ iff $t=a_{i}$ for an $i<n$.

Let 
\begin{eqnarray*}
2+k & := & \max\{ar(f) : f\in\mathcal{F}\}
\\
rk(f) & := & \max\{rk(g)+1 : g<f\}
\\
rk_{2}(f) & := & \max\{rk_{2}(g)+1 : g<f, ar(g)>1\}
\\
m & := & \max\{rk_{2}(f)+1 : f\in\mathcal{F}, ar(f)>1\}
\end{eqnarray*}
Let us denote such a vocabulary $\mathcal{F}$ by $\mathcal{F}^{(m)}_{k}$, and 
the set of terms over $\mathcal{F}^{(m)}_{k}$ by $T^{(m)}_{k}$.
\\

Next let $TIR[\omega^{1+k},\Pi_{2}]$ denote an inference rule
\[
\infer{\forall (x_{0},x_{1},\ldots,x_{k})\, A(x_{0},x_{1},\ldots,x_{k})}
{\forall (x_{0},x_{1},\ldots,x_{k})[\forall (y_{0},y_{1},\ldots,y_{k})<_{lx}(x_{0},x_{1},\ldots,x_{k})\, A(y_{0},y_{1},\ldots,y_{k})\to
A(x_{0},x_{1},\ldots,x_{k})]}
\]
where $A\in\Pi_{2}$ and $(y_{0},y_{1},\ldots,y_{k})<_{lx}(x_{0},x_{1},\ldots,x_{k})$ is the lexicographic ordering
on $(1+k)$-tuples of natural numbers.
Also let $Tr_{\Pi_{1}}$ denote the set of true $\Pi_{1}$-sentences.

Then $TIR^{(m)}_{k}:=TIR^{(m)}[\omega^{1+k},\Pi_{2}]$ [$TIR^{(m)}_{k}+Tr_{\Pi_{1}}:=TIR^{(m)}[\omega^{1+k},\Pi_{2}]+Tr_{\Pi_{1}}$]
denotes a formal system extending the fragment $I\Sigma_{1}$ [$I\Sigma_{1}+Tr_{\Pi_{1}}$], resp.
in which the inference rule $TIR[\omega^{1+k},\Pi_{2}]$ can be applied nestedly at most $m$-times.
To be specific, the derivability relation $TIR^{(m)}_{k}\vdash B$ [$TIR^{(m)}_{k}+Tr_{\Pi_{1}}\vdash B$] is defined recursively as follows.
\begin{enumerate}
\item 
If $I\Sigma_{1}\vdash B$, then $TIR^{(m)}_{k}\vdash B$.

\item
If $I\Sigma_{1}+Tr_{\Pi_{1}}\vdash B$, then $TIR^{(m)}_{k}+Tr_{\Pi_{1}}\vdash B$.

\item
If $TIR^{(m)}_{k}\vdash \forall (x_{0},x_{1},\ldots,x_{k})[\forall (y_{0},y_{1},\ldots,y_{k})<_{lx}(x_{0},x_{1},\ldots,x_{k})\, A(y_{0},y_{1},\ldots,y_{k})\to
A(x_{0},x_{1},\ldots,x_{k})]$, then
$TIR^{(m+1)}_{k}\vdash \forall (x_{0},x_{1},\ldots,x_{k})\, A(x_{0},x_{1},\ldots,x_{k})$.

The same for $TIR^{(m)}_{k}+Tr_{\Pi_{1}}\vdash B$.
\item
If $B$ follows from some $\{C_{i}\}_{i}$ logically and $TIR^{(m)}_{k}\vdash C_{i}$ for any $i$,
then $TIR^{(m)}_{k}\vdash B$.
The same for $TIR^{(m)}_{k}+Tr_{\Pi_{1}}\vdash B$.
\end{enumerate}

\begin{theorem}\label{thm}
For any computable function $f$ the following three conditions are mutually equivalent for each natural number $k$
and each positive integer $m$:
\begin{enumerate}
\item\label{thm1}
$f$ is provably recursive in $TIR^{(m)}_{k}+Tr_{\Pi_{1}}$.

\item\label{thm2}
$f$ is elementary recursive in some fast growing function $F_{\omega^{1+k}\cdot m+q}$ with $q<\omega$.

\item\label{thm3}
There exists a finite term rewriting system $R$ over a vocabulary $\mathcal{F}^{(m)}_{k}$ such that
$R$ reduces under a lpo
with a precedence and
$f$ is elementary recursive in the derivation complexity function
$dh_{R}(n)$.

\end{enumerate}
\end{theorem}

In \cite{attic}, the equivalence of the conditions (\ref{thm1}) and (\ref{thm2}) is shown as Corollary 7.1,
and the implication (\ref{thm2})$\Rightarrow$(\ref{thm3}) is shown in Theorem 8.4.1.
Actually in Definition 8.7 of \cite{attic} the vocabulary $\mathcal{F}^{(m)}_{kQ}=\{list\}\cup\{A_{p}: p<m\}\cup\{f_{q}: q<Q\}$
is introduced with a precedence $list<A_{0}<\cdots<A_{m-1}<f_{0}<\cdots<f_{Q-1}$,
where $list$ is varyadic and $ar(A_{p})=2+k$, $ar(f_{q})=1$.
It is straightforward to see that the implication (\ref{thm2})$\Rightarrow$(\ref{thm3}) holds
when we change the vocabulary to $\mathcal{F}^{(m)}_{kQ}=\{0,S\}\cup\{A_{p}: p<m\}\cup\{f_{q}: q<Q\}$
with the precedence $0<S<A_{p}<f_{q}$ for a constant $0$ and a unary function symbol $S$.

In Theorem 8.3 of \cite{attic} it is shown that for the vocabulary 
$\mathcal{F}^{(m)}_{kQ}=\{list\}\cup\{A_{p}: p<m\}\cup\{f_{q}: q<Q\}$,
$dh_{R}(n)$ is majorized by the slow growing function $G_n(d(\Omega^{2+k}\cdot m+\Omega\cdot Q))$
when $R$ is a finite TRS over $\mathcal{F}_{kQ}^{(m)}$ such that $R$ is reducing under lpo.
Note that $G_n(d(\Omega^{2+k}\cdot m+\Omega\cdot Q))$ is elementary recursive in the fast growing function $F_{\omega^{1+k}\cdot m+Q}$.
Therefore the implication (\ref{thm3})$\Rightarrow$(\ref{thm2}) or equivalently
the implication (\ref{thm3})$\Rightarrow$(\ref{thm1}) is verified partly for the specific vocabulary.

In this note we show the implication (\ref{thm3})$\Rightarrow$(\ref{thm1}) for \textit{any} vocabulary $\mathcal{F}^{(m)}_{k}$
by a proof mining from the proof-theoretic analysis of lpo by Buchholz\cite{Buchholz}.
\\

Let $\mathcal{F}^{(m)}_{k}$ be a finite vocabulary such that the maximal arity of function symbols is $2+k$, and
the maximal rank $rk_{2}(f)=m-1$ for $f\in\mathcal{F}^{(m)}_{k}$.
Also let $R$ be a finite TRS over $\mathcal{F}^{(m)}_{k}$ which is reducing under lpo $<_{lpo}$ 
with a well-founded precedence $<$ on $\mathcal{F}^{(m)}_{k}$.

Let $<_{p}$ be a finite approximation of $<_{lpo}$ defined in p.\,61 of \cite{Buchholz}.
Lemma 7 in \cite{Buchholz} shows that $t\to_{R}s$ is contained in $s<_{p}$ for some $p$.
$p$ depends solely on $R$.
In what follows fix such a $p$.
Let $W_{p}$ denote the well-founded part of $<_{p}$.
There are only finitely many predecessors of any term $t$.
To be specific, the size of the set $\{s: s<_{p}t\}$ is bounded by an elementary recursive function of the size of $t$ (and $p$).
Hence $W_{p}$ is a $\Sigma_{1}$-formula.

\begin{lemma}\label{lem:2}
Let $f\in\mathcal{F}^{(m)}_{k}$ with $n=rk_{2}(f)$, and $n_{2}=n+1$ if $ar(f)>1$.
Otherwise let $n_{2}=n$. Then
$TIR^{(n_{2})}_{k}\vdash\forall a\subset W_{p}(f(a)\in W_{p})$.
\end{lemma}
{\bf Proof}.
By metainduction on $rk(f)$.
We show first
\begin{equation}\label{eq:prg}
TIR^{(n)}_{k}\vdash
\forall a\subset W_{p}[\forall b\subset W_{p}(b<_{lx}a \to f(b)\in W_{p}) \to f(a)\in W_{p}]
\end{equation}
where $b<_{lx}a$ denotes the lexicographic extension of $<_{p}$.

Argue in $TIR^{(n)}_{k}$, and assume that $a\subset W_{p}$ and $\forall b\subset W_{p}(b<_{lx}a \to f(b)\in W_{p})$.
By subsidiary induction on depths of terms $s$ we prove that
$\forall s\in T^{(m)}_{k}(s<_{p}f(a) \to s\in W_{p})$.
This yields $f(a)\in W_{p}$.
Note that the susidiary induction is an instance of $I\Sigma_{1}$.
Let $s<_{p}f(a)$.

If $s\leq_{p}a_{i}$ for an $a_{i}\in a$, then by $a_{i}\in W_{p}$ we have $s\in W_{p}$.

Next consider the case $s=f(b)$ with $b<_{lx}a$.
By SIH and $a\subset W_{p}$, we have $b\subset W_{p}$.
The assumption yields $s\in W_{p}$.

Finally consider the case $s=g(b)$ with $g<f$.
Then $rk_{2}(g)\leq rk_{2}(f)$, and $rk_{2}(g)=rk_{2}(f) \Rightarrow ar(g)\leq 1$.
We have $TIR^{(n)}_{k}\vdash\forall b\subset W_{p}(g(b)\in W_{p})$ by MIH.
On the other hand we have $b\subset W_{p}$ by SIH.
Hence $g(b)\in W_{p}$.
Thus (\ref{eq:prg}) is shown.

If $ar(f)\leq 1$, then (\ref{eq:prg}) yields $TIR^{(n)}_{k}\vdash\forall s\in W_{p}(f(s)\in W_{p})$.
Let $2+\ell:=ar(f)>1$ with $\ell\leq k$.
Let $G(a):\Leftrightarrow \forall t\in W_{p}(f(a*(t))\in W_{p})$,
where $a$ ranges over sequences of terms in length $1+\ell$ and
$(t_{0},t_{1},\ldots,t_{\ell})*(t)=(t_{0},t_{1},\ldots,t_{\ell},t)$.
$G(a)$ is a $\Pi_{2}$-formula.
Next we show
\begin{equation}\label{eq:prg2}
TIR^{(n)}_{k}\vdash \forall a\subset W_{p}[\forall b\subset W_{p}(b<_{lx}a \to G(b)) \to G(a)]
\end{equation}
Then an application of $TIR[\omega^{1+\ell},\Pi_{2}]$ yields
$TIR^{(n+1)}_{k}\vdash \forall a\subset W_{p}\, G(a)$, i.e., 
$TIR^{(n+1)}_{k}\vdash\forall t_{0},\ldots,t_{\ell+1}\in W_{p}(f(t_{0},\ldots,t_{\ell+1})\in W_{p})$ as desired.

Argue in $TIR^{(n)}_{k}$.
Assume $a\subset W_{p}$ and $\forall b\subset W_{p}(b<_{lx}a \to G(b))$.
We need to show $G(a)$, i.e., $\forall t\in W_{p}(f(a*(t))\in W_{p})$.
We show this by induction on $t\in W_{p}$.
Again this is an instance of $I\Sigma_{1}$.
Let $t\in W_{p}$, and suppose that
$\forall s<_{p}t(f(a*(s)\in W_{p})$ as IH.
We need to show $f(a*(t))\in W_{p}$.
By (\ref{eq:prg}) it suffices to show that
$\forall b\subset W_{p}\forall s\in W_{p}(b*(s)<_{lx}a*(t) \to f(b*(s))\in W_{p})$.
Let $b\subset W_{p}$, $s\in W_{p}$ and $b*(s)<_{lx}a*(t)$.
If $b<_{lx}a$, then the assumption yields $G(b)$ and $s\in W_{p}$.
Hence $f(b*(s))\in W_{p}$.
Let $b=a$. Then $s<_{p}t$.
IH yields $f(b*(s))\in W_{p}$.

This shows (\ref{eq:prg2}), and a proof of Lemma \ref{lem:2} is completed.
\hspace*{\fill} $\Box$
\\

\begin{lemma}\label{lem:3}
$TIR^{(m)}_{k}\vdash\forall t\in T^{(m)}_{k}(t\in W_{p})$.
\end{lemma}
{\bf Proof}.
By induction on depths of terms $t\in T^{(m)}_{k}$ using Lemma \ref{lem:2}.
Note that this is an instance of an axiom $I\Sigma_{1}$.
\hspace*{\fill} $\Box$
\\

Lemma \ref{lem:3} shows that the derivation complexity function $dh_{R}(n)$ is provably recursive in $TIR^{(m)}_{k}$.
This yields the implication (\ref{thm3})$\Rightarrow$(\ref{thm1}) in Theorem \ref{thm}.

\end{document}